\documentstyle[12pt]{article} 

\def\xxx#1{}	\def\mat{\def}	\def\grk{\def}	\def\pol{\def}
	
\xxx{
\setlength{\textwidth}{164mm}	\setlength{\oddsidemargin}{2mm}
\setlength{\textheight}{230mm}	\addtolength{\topmargin}{-2cm}
}

\mat\0#1{\makebox{\(#1\)}} 
\mat\8#1{\makebox{\(\displaystyle #1\)}} 

\mat\9#1{\0{ \left( {#1} \right) }}
\mat\1#1{\0{ \left| {#1} \right|  }}	
\mat\2#1{\0{ \left\|{#1} \right\| }}
\mat\3#1{\0{ \left\langle {#1} \right\rangle}}	
\mat\4#1{\0{ \left[ {#1} \right]  }}	
\mat\5#1{\0{ \left\{{#1} \right\} }}	
\mat\6#1#2{\0{#1\in #2}}
\mat\7#1#2#3{\0{#1:#2\to #3}}	

\xxx{
\pol\cu{\raisebox{-1.45ex}[-2ex][-2ex]{``}}
\polae#1#2#3{\makebox[0em][l]{\hspace*{#1}\raisebox{#2}[0mm][0mm]{`}}#3}
\pol\a.{\polae{.31em}{-1.45ex}{a}}
\pol\A.{\polae{.44em}{-1.45ex}{A}}
\pol\e.{\polae{.20em}{-1.45ex}{e}}
\pol\E.{\polae{.41em}{-1.45ex}{E}}
\pol\u.{{\l}}	
} 

\mat\AA{{\cal A}}	\mat\JJ{{\cal J}}	\mat\SS{{\cal S}}	
\mat\BB{{\cal B}}	\mat\KK{{\cal K}}	\mat\TT{{\cal T}}	
\mat\CC{{\cal C}}	\mat\LL{{\cal L}}	\mat\UU{{\cal U}}	
\mat\DD{{\cal D}}	\mat\MM{{\cal M}}	\mat\VV{{\cal V}}	
\mat\EE{{\cal E}}	\mat\NN{{\cal N}}	\mat\WW{{\cal W}}	
\mat\FF{{\cal F}}	\mat\OO{{\cal O}}	\mat\XX{{\cal X}}	
\mat\GG{{\cal G}}	\mat\PP{{\cal P}}	\mat\YY{{\cal Y}}	
\mat\HH{{\cal H}}	\mat\QQ{{\cal Q}}	\mat\ZZ{{\cal Z}}	
\mat\II{{\cal I}}	\mat\RR{{\cal R}} 

\mat\N{{\bf N}}		\mat\Z{{\bf Z}}		
\mat\R{{\bf R}}		\mat\C{{\bf C}}

\grk\Gam{\Gamma}	\grk\Del{\Delta}	\grk\Th{\Theta}   
\grk\Lam{\Lambda}       \grk\Sig{\Sigma}	\grk\Ups{\Upsilon}	
\grk\Yps{\Upsilon}    	\grk\Ph{Phi}		\grk\Ps{\Psi} 
\grk\Om{\Omega}     

\grk\al{\alpha}	 	\grk\gam{\gamma}	\grk\del{\delta}	
\grk\eps{\varepsilon} 	\grk\th{\vartheta}	\grk\kap{\kappa}	
\grk\lam{\lambda}	\grk\rh{\varrho}	\grk\sig{\sigma}	
\grk\ups{\upsilon} 	\grk\ph{\varphi}	\grk\ps{\psi} 		
\grk\om{\omega}

\xxx{
\grk\alfa{\alpha} 	\grk\dzeta{\zeta} 	\grk\dz{\zeta} 
\grk\teta{\vartheta}	\grk\jota{\iota} 	\grk\mi{\mu}	
\xxx{\grk\ni{\nu}} 	
\grk\ksi{\xi} 		\grk\ro{\varrho} 	
\grk\yps{\ypsilon}
} 

\newtheorem{xDef}{Definition}

\newtheorem{xTheo}{Theorem}
\newtheorem{xCor}{Corollary}
\newtheorem{xLem}{Lemma}
\newtheorem{xProp}{Proposition}
\newenvironment{Def}{\begin{xDef}\rm }{\end{xDef}}
\newenvironment{Theo}{\begin{xTheo}\sl }{\end{xTheo}}

\newenvironment{Lem}{\begin{xLem}\sl }{\end{xLem}}
\newenvironment{Prop}{\begin{xProp}\sl }{\end{xProp}}

\mat\Pf{{\em Proof. }}       	
\mat\fP{\rule{2.5mm}{2.5mm}\vspace*{1.5ex}}
\mat\barr{\begin{array}} 	\mat\earr{\end{array}} 
\mat\bequ{\begin{equation}} 	\mat\eequ{\end{equation}} 
\mat\bequarr{\begin{eqnarray}}  \mat\eequarr{\end{eqnarray}}
\mat\bite{\begin{itemize}} 	\mat\eite{\end{itemize}}

\mat\Then{\Longrightarrow}
\mat\to{\rightarrow}		\mat\from{\leftarrow}
\mat\To{\longrightarrow}	\mat\From{\longleftarrow}
\mat\onto{\stackrel{{\scriptsize\rm onto}}{\to}}

\mat\stk{\stackrel}		\mat\fr#1#2{{\textstyle\frac{#1}{#2}}}
\mat\xl#1{\mbox{\rm\ \ #1\ \ }}
\mat\r#1{\mbox{\rm #1}}
\mat\til#1{\0{\widetilde{#1}}}	\mat\ovl#1{\0{\overline{#1}}}

\mat\d.{\partial}	
\mat\c.{\cdot}		\mat\x.{\times}	

\mat\eq{\equiv}		\mat\app{\approx} 
\mat\ctin{\subset}      \mat\cts{\supset}
\mat\mns{\setminus} 	
\mat\Union{\Bigcup}     \mat\Inters{\Bigcap}

\mat\dg{\mbox{deg$_2$\,}}
\mat\dist{\mbox{dist\,}} 
\mat\cl{\mbox{cl\,}} 
\mat\pitchfork{\top}
\mat\e.{\epsilon}

\begin{document}\sloppy

\title{Cyclic hybrid systems of flows on manifolds    
\footnote{Supported by Polish MNiSzW Grant No.\ 1\,P03A\,015\,29.}
}
\author{Witold Szczechla} 
\date{} 
\maketitle

\begin{abstract} The considered continuous-and-discrete hybrid system is a 
cyclic relay of smooth flows on an $n$-dimensional manifold $M$, where the 
discrete process of switching from each flow to the next takes place on 
the boundaries of some fixed compact $n$-dimensional submanifolds of $M$. 
The main result is the existence of at least one periodic cyclic 
trajectory under some topological condition concerning one of the 
submanifold boundaries.  \end{abstract}

\section{Introduction} 

Whereas single vector fields on manifolds and their local combinations 
have been extensively studied, especially in the areas of theoretical 
mechanics, dynamical systems and differential geometry,  
the complexity of the {\em global behaviour of combined trajectories\/} 
has attracted less research, even in control-theoretic setting.     
Naturally arising in this context are relays of flows on manifolds,   
defined as special smooth-and-discrete hybrid systems 
in \cite{Dis}. They constitute a special case of 
nonlinear systems studied, in particular, in \cite{S1}, \cite{Dis}, 
\cite{ProcHS}, \cite{Three}, \cite{2Nper}.  

A {\em relay of flows\/} is a system whose phase space is a smooth 
($C^\infty$) manifold $M$ of dimension $n>0$. Let $p$ be an integer, $p>1$. 
One is given a cyclic collection of smooth flows on $M$, 
	\[ \7{F_k}{M\times\R}{M} \ \ \ \ \9{k\in\Z_p}, 
	\]   
and a corresponding collection of compact $n$-dimensional 
manifolds $M_k$ with nonempty boundaries \9{k\in\Z_p}. The following 
two conditions are assumed:   
	\begin{itemize}  
	\item[{(i)}]  For each $k$ we have 
\[   \d.M_{k-1}\cap M_k=\emptyset.  
\]  
	\item[{(ii)}]  For each $k$ there is a value $T_k>0$ such that  
\[  F_k(T_k,\d.M_{k-1})\ctin\r{int}M_k. 
\] 
	\end{itemize} 

The cyclic hybrid system thus obtained will be called $\cal S$. A {\em 
trajectory\/} (i.e., a {\em solution\/} or more precisely, a   
{\em quasisolution\/}) of the system $\cal S$ is defined to be a pair of 
functions \9{F,\alpha} with \7{F}{[t_0,\infty)}{M} continuous and 
\7{\alpha}{[t_0,\infty)}{\Z_p} piecewise constant. The following 
two conditions must be satisfied:
	\begin{itemize}  
	\item[{(i)}]  If $\al|[a,b]=k=\r{const}$, then 
		\[ F(b) = F_k^{b-a}(F(a)). \] 
	\item[{(ii)}]  If $\al(t-0)=k\neq\al(t+0)$, then $\al(t+0)=k+1$ 
and $F(t)\in\d.M_k$.  
	\end{itemize} 

{\em Note.\/} One speaks of `quasisolutions' since the term `solution' is 
sometimes reserved for the trajectories satisfying the following 
additional condition, not considered in this paper: (iii) if $\al(t)=k$, 
then $F(t)\not\in\r{int}M_k$.

{\em Motivation.\/} Besides its theoretical interest, the model has been 
motivated by the switching occurring in electromagnetc fields and/or some 
designs of control systems engineering.


\subsection{The limit sets for solutions} 

For any initial point $x_0\in M$ and initial mode $\al(0+)=k_0$ we define 
the accessible set $A(x_0,k_0)$ and the $\om$-limit set $\Om(x_0,k_0)$ for 
{\em solutions\/} of $\cal S$ in the usual way (by analogy to dynamical 
systems, cf. \cite{K-H}).  We note that trajectories can be nonunique and 
any point in the $\om$-limit set has (by definition) to be approached by 
{\em some\/} trajectory starting at the end of {\em every\/} trajectory of 
$x_0$ which has just completed $m$ switches (where $m$ is arbitrary).

\begin{Prop}  Each of the sets $A(x_0,k_0)$ and $\Om(x_0,k_0)$ is 
connected.  
\end{Prop}

\Pf The set $A(x_0,k_0)$ is made up of consecutive pieces of trajectories 
between $\d.M_{k-1}$ and $\d.M_k$ which may be branching, but connect at 
the ends.  Hence, the closure of $A(x_0,k_0)$ is always connected.  By 
definition, the limit set $\Om(x_0,k_0)$ can be expressed as the 
intersection of a descending sequence of such closures.  \fP


\section{Periodic quasisolutions}   

In general, the system $\cal S$ may have no periodic quasisolutiotion.  
Such may be the case if $p=2$, $M=S^2$, and the $M_k$ are two disjoint 
disks whose boundaries are twisted by the first return map.  In such 
a case, however, the disk $M_0$ can never be entirely carried inside $M_1$ 
by the flow $F_1$ because of an argument using Brouwer fixed point 
theorem.  In fact, the following specific property of the sphere $S^{n-1}$, 
related to the degree (mod~2) of a mapping, 
turns out to be relevant to our method. 
	\begin{Def} Let $S$ be a compact manifold.  We will say that $S$ 
has the {\em coincidence property (mod~2)\/} if $S$ is connected and if 
for every compact manifold $N$ of the same dimension, for every pair of smooth 
maps \7{f,g}{N}{S} satisfying the condition $\dg(f)\neq\dg(g)$,  
there exists a point of coincidence 
(i.e., a point $x\in N$ such that $f(x)=g(x)$).
	\end{Def} 
Let us notice that the sphere of any dimension,
$S=S^{n-1}$, has the above property.  Indeed, suppose $f(x)\neq 
g(x)$ for all $x\in N$.  Then $f$ is homotopic to the map $-g=\sigma\circ 
g$, where $\sigma$ is the antipodal map.  However, since $\dg(\sigma)=1$, 
this implies $\dg(f)=\dg(-g)=1\cdot\dg(g)$.

\begin{Theo} \label{relay1}  Suppose that, for some index $\beta\in\Z_p$,  
we have 
\bequ\label{Mbeta} 
F_\beta(t,M_{\beta-1}) \ctin \r{int}M_\beta \xl{for all}\  t\geq T_\beta.  
\eequ 
Suppose also that, for some index $k\in\Z_p$, the boundary $\d.M_k$ has 
the coincidence property (mod~2). 
Then the system $\cal S$ has a $p$-periodic quasisolution.
\end{Theo}   

We may, and will, assume that the manifold $\d.M_0$ has the coincidence 
property (mod~2). Let $x\in\d.M_{k-1}$, where $k$ is arbitrary.  By the 
assumptions (i) and (ii), it is reasonable to suppose that, generically, 
there are an odd number of values $t\in(0,T_k)$ for which 
$F_k(t,x)\in\d.M_k$, since $x\not\in M_k$ and $F_k(T_k,x)\in\r{int}M_k$. 
On the other hand, by assumption~(\ref{Mbeta}), the {\em backward\/} 
trajectory of the flow $F_\beta$ starting at a generic point 
$x\in\d.M_\beta$ should cross $\d.M_{\beta-1}$ an even number of times as 
neither $x$ nor $F_\beta(-T_k,x)$ can belong to $M_{\beta-1}$. We want to 
relate the above observations to the degree mod~2 of certain maps into 
$\d.M_0.$


\subsection{Technical lemmas}

For technical reasons, we will approximate our system of submanifolds by a 
certain family of similar systems continuously parametrized by 
$\lam\in\R^{p+1}$ and prove a suitably generalized theorem for a residual 
set of parameters. The new system will consist of manifolds $M_i^\lam\ \ 
(i=0,\ldots p),$ where $M_i^\lam$ is close to $M_i,$ but possibly 
$M_0^\lam\not= M_p^\lam.$ In particular, $i$ is not a cyclic index here. 
The flows will be unchanged and indexed $1$ to $p$.

For $i=0,1,\ldots,p,$ pick a smooth function $f_i:M\to\R$ having the
following properties: 

\begin{itemize} 
\item[(1)] $f_i^{-1}(\R_{+})=M_i$ 
\item[(2)] 0 is a regular value of $f_i.$
\item[(3)] If $\lim_{n\to\infty}x_n=0$ then $\lim_{n\to\infty}\dist 
(x_n,\d.M_i)=0.$
\end{itemize}
We note that, by (2) and (3), every value sufficiently close to 0 is
regular.

\begin{Lem}
There exist functions $f_i$ satisfying (1)--(3).
\end{Lem}
\Pf
First we take flows $\Psi_i,$ transverse to $\d.M_i$ and directed
outwards. Next, define the  $f_i$ on neighborhoods of $\d.M_i$ by 
\bequ\label{3151}
f_i\90{\Psi_i(t,\xi)}=-\nu(t),\ \ \ (\xi\in\d.M_i,\ |t|<\e.), 
\eequ
where $\e.$ is so small that the maps $(t,\xi)\mapsto\Psi_i(t,\xi)$ are
diffeomorphic embeddings of $(-\e.,\e.)\times\d.M_i$ into $M$ and 
\[
\nu(t)=\left\{
\begin{array}{lll}
t & {\rm for}& |t|<\e./4\\ 
\e./3 & {\rm for} & |t|>\e./2.
\end{array}
\right.
\]
Now it is enough to extend each $f_i$ by $\e./3$ on $M_i$ and $-\e./3$
on $M\mns M_i.$\fP

For $\lam=(\lam_0,\lam_1,\ldots,\lam_p)\in\R^{p+1},$ denote
$f_i^{-1}\90{[\lam_i,\infty)}$ by $M_i^\lam.$ Let $\Lam\ctin\R^{p+1}$ be
the set of those $\lam=(\lam_0,\lam_1,\ldots,\lam_p)$ which satisfy the
following conditions: 
	\begin{itemize} 
	\item[(4)] For each $i,$ $\lam_i$
is a regular value of $f_i$ 
	\item[(5)] The sets $M_i^\lam$ satisfy the
hypothesis of our theorem with $M_i^\lam$ playing the role of $M_i$ 
except possibly for the requirement $M_0^\lam=M_p^\lam.$ 
	\item[(6)]
$F_i(T_i,\d.M_{i-1}^\lam)\ctin\r{int}M_i^\lam\ (i=1,\ldots,p)$ and
$F_{\beta}(T_{\beta},M_{\beta-1}^\lam)\ctin\r{int}M_{\beta}^\lam.$ 
\end{itemize}

Using properties (1)--(3), it is easy to show that $\Lam$ is an open
set and that $0\in\cl\Lam.$ For the remainder of the proof it will be
usually understood that $\lam=(\lam_0,\lam_1,\ldots,\lam_p)\in\Lam.$ A
generic property means a property holding on a residual subset of
$\Lam.$


\subsection{The manifold of the switching points} 

Let us define the set $N^{\lam}\ctin M\times\R^p$ to be the set of the 
points $(x,t_1,\ldots,t_p)$ satisfying: 
  \begin{itemize}
    \item[(7)] $0<t_i<T_i$ for each $i$;
    \item[(8)] if we denote $x_0=x$ and $x_{i+1}=F_{i+1}(t_{i+1},x_i)$, 
then $x_i\in\d.M_i^\lam$ for each $i.$
  \end{itemize}
Notice that $N^\lam$ is always nonempty and compact.  (While we do 
admit the possibility of $N^\lam$ being disconnected, the essential 
assumption is the connectivity of $\d.M_0$.)  

	\begin{Lem}\label{Nlambda} Generically, $N^\lam$ is a smooth 
$(n-1)$-dimensional submanifold of $M\times\R^p.$ 
	\end{Lem} \Pf Let
$\Omega=M\times(0,T_1)\times\ldots\times(0,T_p).$ Define the map 
$\nu:\Omega\to\R^{p+1}$ in the following way. Let 
$\omega=(x,t_1,\ldots,t_p)\in\Omega,$ $x_0=x$ and 
$x_{i+1}=F_{i+1}(t_{i+1},x_i)\ 
(i=0,\ldots,p-1);$ then \[ 
\nu(\omega)=\90{f_0(x_0),f_1(x_1),\ldots,f_p(x_p)}.  \] Then 
$N^\lam=\nu^{-1}(\lam),$ so $N^\lam$ is a smooth submanifold of dimension 
$\dim\Omega-\dim\R^{p+1}=n-1$ whenever $\lam$ is a regular value of $\nu.$ 
By Morse-Sard theorem, the set of the regular values is residual in 
$\R^{p+1}$ and hence in $\Lam.$ \fP


\subsection{Computation of the degrees mod 2}

Define the maps $\nu_0:N^\lam\to\d.M_0^\lam$ and 
$\nu_1:N^\lam\to\d.M_p^\lam$ 
by
\[ \nu_0(x,t_1\ldots t_p)=x  \]
and
\[ \nu_1(x,t_1\ldots t_p)=F^{t_p}_{p}\circ\cdots\circ
F^{t_1}_1(x), \] 
where we used the notation $F_i^t(z)=F_i(t,z)$.
  
Notice that Lemma~\ref{Nlambda}, together with the connectivity of of the 
manifolds  
$\d.M_0$ and $\d.M_p^\lam$, allow us to compute, for almost every 
$\lam$, the degree modulo 2 of $\nu_0$ and $\nu_1$.      
For any $x\in M$ we put $R_i^i(x)=\{x\}$ and define the sets
$R_i^j(x)\ \ (i,j=0,\ldots,p)$  
recursively by
\[  
R_i^{j+1}(x)=\bigcup_{0<t<T_{j+1}}F^t_{j+1}\90{R_i^j(x)}\cap\d.M_{j+1}^\lam\ 
\ 
        {\rm for}\ \ i\leq j<p \] and
\[  
R_i^{j-1}(x)=\bigcup_{0<t<T_j}F^{-t}_{j}\90{R_i^j(x)}\cap\d.M_{j-1}^\lam
     \ \ {\rm for}\ \ 0 < j\leq i.  \] 
Also, let $V_i$ denote the vector field inducing the flow $F_i.$ 


\begin{Lem}\label{1-0} Generically, if $N^\lam$ is a smooth 
$(n-1)$-dimensional submanifold, then $\dg\nu_0=1$ and $\dg\nu_1=0.$ 
\end{Lem}
\Pf Let $\Lam_0$ be the set of those $\lam$ for which there exists an
$x\in\d.M_0^\lam$ with the following transversality property:
\begin{equation}\label{9} {\rm If}\ \ y\in R_0^j(x)\ \ {\rm and}\ \
j>0\ \ {\rm then}\ \ V_j(y)\pitchfork\d.M_j^\lam.  \end{equation}
Since $F_i(T_i,\d.M_{i-1}^\lam)\cap\d.M_i^\lam=\emptyset,$ it is
clear that the set of pairs $(\lam,x)$ satisfying (\ref{9}) is open. In
particular, the set $\Lam_0$ is open.

We will show that $\Lam_0$ is dense. Take a point $z\in\d.M_0^\lam,$ 
choose
$\xi_1$ close to $\lam_1$ as a regular value of the map 
$g(t)=f_1\90{F_1(t,z)}$ 
and
put  $S_1=\{F_1(t,z)\mid t\in g^{-1}(\xi_1)\}.$
Having defined $\xi_i$ and $S_i$ (where $i<p),$ choose
$\xi_{i+1}$ close to $\lam_{i+1}$ as a common regular value of the maps
  \[ g_y(t)=f_{i+1}\90{F_{i+1}(t,y)}  \]
for all $y\in S_i,$ and put 
  \[ S_{i+1}=\{F_{i+1}(t,y)\mid y\in S_i, t\in g_y^{-1}(\xi_{i+1})\}. \]
(Note that each $S_i$ will be at most countable.) Condition~(\ref{9}) will be 
satisfied by the substitution $x=z$ and $\lam=(\lam_0,\xi_1\ldots\xi_p).$ 

Suppose $\lam\in\Lam_0$ and take an $x\in\d.M_0^\lam$ 
satisfying~(\ref{9}). By Morse-Sard theorem, we may take $x$ to be a 
regular value of $\nu_0.$ The number $\dg\nu_0$ is equal to the reduction 
modulo 2 of the number of points in $\nu_0^{-1}(x).$ Consider any $y\in 
R_0^j(x),$ where $0\leq j<p.$ Since $y\notin M_{j+1}^\lam$ and 
$F_{j+1}(T_{j+1},y)\in\r{int} M_{j+1}^\lam,$ the transversality condition 
implies that the set \[ \{ t\in(0,T_{j+1})\mid 
F_{j+1}(t,y)\in\d.M_{j+1}^\lam \} \] has an odd number of elements, which 
proves $\nu_0^{-1}(x)$ has an odd number of points, so $\dg\nu_0=1.$

Let $\Lam_1$ be the set of those $\lam$ for which there exists an
$x\in\d.M_p^\lam$ with the following property:
\begin{equation}\label{10} {\rm If}\ y\in R_p^j(x)\ {\rm and}\ j<p\
{\rm then}\ V_j(y)\pitchfork \d.M_j^\lam.  \end{equation} As before,
we first observe that $\Lam_1$ is open and dense. Now suppose $\lam$ and
$x\in\d.M_p^\lam$ satisfy~(\ref{10}) and $x$ is a regular value of
$\nu_1.$ Consider any $y\in R_p^{\beta}(x).$ Since $y\in
M^\lam_{\beta},$ $y\notin M_{\beta-1}^\lam.$ Since $y\notin\r{int} 
M_{\beta}^\lam$ and 
$F^{\beta}(T_{\beta},\beta)(M_{\beta-1}^\lam) 
\ctin \r{int}M_\beta^\lam,$ 
also
$F^{-T_{\beta}}_{\beta}(y)\notin M_{\beta-1}^\lam.$ Thus, the
transversality condition implies that the set 
\[ 
\{ t\in(0,T_{\beta})\mid F_{\beta}(-t,y)\in\d.M_{\beta-1} \} 
\] 
has an even number of elements.  Now consider any $y\in R_p^{j+1}(x),$
where $\beta-1>j>0.$ Since $\d.M_{j+1}^\lam$ is connected and
$F_{j+1}(-T_{j+1},\d.M_{j+1}^\lam)\cap\d.M_j^\lam=\emptyset,$ the
set $F_{j+1}(-T_{j+1},\d.M_{j+1}^\lam)$ must be contained in 
$\r{int}M_j^\lam$ or $M\mns M_j^\lam.$ Similarly, $\d.M_{j+1}^\lam$ is 
contained
in $\r{int}M_j^\lam$ or $M\mns M_j^\lam,$ since $\d.M_{j+1}^\lam\cap\d.
M_j^\lam=\emptyset.$ The transversality condition then implies that the
set 
\[ 
\{ t\in(0,T_{j+1})\mid F_{j+1}(-t,y)\in\d.M_j^\lam \} 
\] 
has an even or odd number of elements independently of $y.$ Thus we
have proved that $\nu_1^{-1}(x)$ has an even number of points, which
means that $\dg\nu_1=0.$ \fP


Now, for any $\lam\approx 0$ there is a diffeomorphic projection    
\7{\pi_\lam}{\d.M_p^\lambda}{\d.M_0} close to the identity.  
It follows that the maps \7{f,g}{N^\lam}{\d.M_0}, where $f=\nu_0$ 
and $g=\pi_\lam\circ\nu_1$ still satisfy $\dg(f)\neq\dg(g)$.  The assumed 
coincidence property (mod~2) of the manifold $\d.M_0$ implies that 
$f(x_\lam)=g(x_\lam)$ for some point $x_\lam\in N^\lam$.  Thus, for 
each value of $\lam\approx 0$ there exists a trajectory which is nearly 
closed (up to the map $\pi_\lam$).  By letting $\lam\to 0$ and selecting a 
suitable convergent subsequence, we obtain a $p$-periodic trajectory of 
the hybrid system $\cal S$.  This completes the proof of 
Theorem~\ref{relay1}.  \fP

\vspace*{3ex} 

\parbox{8cm}{{\sc Instytut Matematyki\\ Uniwersytet Warszawski\\ 
ul.\ Banacha 2\\ 00-913 Warszawa, Poland}\\ 
E-mail address: {\tt witold@mimuw.edu.pl}}

\end{document}